\documentclass[12pt,reqno]{amsart}
\usepackage{verbatim}
\usepackage{amscd, amsfonts}

\begin{document}

\numberwithin{equation}{section}
\renewcommand{\theequation}{\thesection.\arabic{equation}}
\setcounter{secnumdepth}{2}

% The following establish some abbreviations for combinations and/or
% commonly used letters/terms.
%
% The following are macros for Greek letters & special fonts.
\newcommand{\Ac}{{\mathcal A}}
\newcommand{\Bc}{{\mathcal B}}
\newcommand{\bt}{\tilde{b}}

\newcommand{\bfA}{\bf{A}}
\newcommand{\bfB}{\bf{B}}
\newcommand{\bfF}{\bf{F}}
\newcommand{\BH}{\Bc_H}

\newcommand{\Cc}{{\mathcal C}}
\newcommand{\Dc}{{\mathcal D}}
\newcommand{\Ec}{{\mathcal E}}

\newcommand{\Fc}{{\mathcal F}}

\newcommand{\Gc}{{\mathcal G}}
\newcommand{\GD}{\Gc_D}
\newcommand{\Gma}{{\Gamma}}

\newcommand{\hbf}{\bf{h}}
\newcommand{\Hc}{{\mathcal H}}
\newcommand{\htil}{\tilde{h}}
\newcommand{\kbf}{\bf{k}}
\newcommand{\Lap}{\Delta}
\newcommand{\Mcl}{{\mathcal M}}

\newcommand{\nab}{\nabla}
\newcommand{\Nc}{{\mathcal N}}
\newcommand{\Oc}{\mathcal{O}}
\newcommand{\Om}{\Omega}

\newcommand{\pal}{\partial}
\newcommand{\Pc}{\mathbb{P}}

\newcommand{\Rb}{\overline{\R}}
\newcommand{\R}{\mbox{$\mathbb R$}}

\newcommand{\Sc}{{\mathcal S}}
\newcommand{\sg}{\sigma}
\newcommand{\sgt}{\tilde{s}}
\newcommand{\vap}{\varphi}

\newcommand{\vbf}{\bf{v}}
\newcommand{\wbf}{\bf{w}}

\newcommand{\Wc}{{\mathcal W}}

\newcommand{\barry}{\begin{eqnarray}}
\newcommand{\bc}{\begin{center}}
\newcommand{\beq}{\begin{equation}}
\newcommand{\bpf}{\begin{proof} \quad}
\newcommand{\btm}{\begin{thm} \quad}

\newcommand{\earry}{\end{eqnarray}}
\newcommand{\ec}{\end{center}}
\newcommand{\eeq}{\end{equation}}
\newcommand{\epf}{\end{proof}}
\newcommand{\etm}{\end{thm}}

\newcommand{\deq}{:= }
\newcommand{\deqs}{\ :=\ }

\newcommand{\eqs}{\ =\ }
\newcommand{\geqs}{\ \geq \ }
\newcommand{\leqs}{\ \leq \ }
\newcommand{\mns}{\, - \, }
\newcommand{\opluss}{\ \oplus  \ }
\newcommand{\pls}{\, + \, }
\newcommand{\plms}{\ + \ }

% The following are standard combinations of letters.

\newcommand{\foral}{\qquad \mbox{for all} \quad }
\newcommand{\foreach}{\qquad \mbox{for each} \quad }

\newcommand{\Pslemma}{\ \mbox{Poincar\'{e}'s lemma} \ }
\newcommand{\wrt}{ \mbox{with respect to}  }
\newcommand{\xand}{\quad \mbox{and} \quad }
\newcommand{\xfor}{ \quad \mbox{for} \ }
\newcommand{\xiff}{\ \mbox{if and only if} \ }

\newcommand{\xon}{\qquad \mbox{on} \quad}
\newcommand{\xor}{\qquad \mbox{or} \quad}
\newcommand{\xthen}{\quad \mbox{then} \quad}
\newcommand{\xwhen}{\qquad \mbox{when} \quad}
\newcommand{\xwith}{\qquad \mbox{with} \quad}

% The following macros are for regions and other subsets of R^N.  Line 130
\newcommand{\bdy}{\partial \Omega}

\newcommand{\Gamo}{\Gamma_1}
\newcommand{\Gamj}{\Gamma_j}
\newcommand{\GamJ}{\Gamma_J}
\newcommand{\Gamk}{\Gamma_k}
\newcommand{\Gamx}{\Gamma_x}

\newcommand{\gamo}{\gamma_1}
\newcommand{\gamj}{\gamma_j}
\newcommand{\gamk}{\gamma_k}
\newcommand{\gamL}{\gamma_L}

\newcommand{\Gamz}{\Gamma_0}
\newcommand{\Gamnu}{\Gamma_{\nu}}
\newcommand{\Gamtau}{\Gamma_{\tau}}
\newcommand{\Gamtil}{\tilde{\Gamma}}

\newcommand{\Omt}{\Om \times \Om}
\newcommand{\Omext}{\Om_e}

\newcommand{\Rn}{{\R}^n}
\newcommand{\RN}{{\R}^N} 
\newcommand{\RNp}{{\R}^N_+} 
\newcommand{\Rte}{{\R}^3}

\newcommand{\Rto}{{\R}^2}
 
 % The following macros are for measures and integrals.
\newcommand{\dsg}{\, ds}
\newcommand{\dtx}{\, d^2x}

\newcommand{\Iby}{\int_{\bdy} \; }
\newcommand{\IGamnu}{\int_{\Gamnu}}
\newcommand{\IGamtau}{\int_{\Gamtau}}
\newcommand{\IGamt}{\int_{\Gamtil}}

\newcommand{\IOm}{\int_{\Om} \; }

% The following macros are for linear or bilinear  operators and maps.  This is line 144
\newcommand{\ang}[1]{\langle#1\rangle}
\newcommand{\curlphi}{\gradp \vap}
\newcommand{\curlpsi}{\gradp \psi}

\newcommand{\divF}{\Div \bfF}
\newcommand{\divG}{\Div \bfG}

\newcommand{\gamu}{\gamma(u)}
\newcommand{\gamphi}{\gamma(\vap)}
\newcommand{\gampsi}{\gamma(\psi)}

\newcommand{\gradp}{\nabla^{\perp}}

\newcommand{\gradchi}{\nabla \chi}
\newcommand{\gradphi}{\nabla \varphi}
\newcommand{\gradphiz}{\nabla \vapz}
\newcommand{\gradpsi}{\nabla \psi}
\newcommand{\gradpchi}{\gradp \chi}

\newcommand{\gradxi}{\nabla \xi}
\newcommand{\gradu}{\nabla u}
\newcommand{\gradv}{\nabla v}

\newcommand{\Lapbj}{\Delta \, b_j}
\newcommand{\Lapchi}{\Delta \chi}
\newcommand{\Lapphi}{\Delta \vap}
\newcommand{\Lappsi}{\Delta \psi}

\newcommand{\Lapu}{\Delta u}
\newcommand{\Lapv}{\Delta v}

\newcommand{\curl}{\mathop{\rm curl}\nolimits}
\newcommand{\Curl}{\rm Curl}

\newcommand{\Div}{\mathop{\rm div}\nolimits}
\newcommand{\Dnu}{D_{\nu}}

\newcommand{\PG}{P_G}
\newcommand{\PGz}{P_{G0}}
\newcommand{\PC}{P_C}
\newcommand{\PCz}{P_{C0}}
\newcommand{\QG}{Q_G}
\newcommand{\QGz}{Q_{G0}}
\newcommand{\QC}{Q_C}
\newcommand{\QCz}{Q_{C0}}

% The following macros are  for scalar function spaces.   line 166
\newcommand{\COm}{C(\Om)}
\newcommand{\Ccone}{C_c^1(\Om)}
\newcommand{\COmz}{C_0(\Om)}

\newcommand{\HGamnu}{H^1_{\Gamnu 0}(\Om)}
\newcommand{\HGamtau}{H^1_{\Gamtau 0}(\Om)}
\newcommand{\HGamz}{H^1_{\Gamma 0}(\Om)}

\newcommand{\Hone}{H^1(\Om)}
\newcommand{\Hmone}{H^1_m(\Om)}
\newcommand{\Hnone}{H^{-1}(\Om)}

\newcommand{\Hstar}{\Hone^*}

\newcommand{\Hzone}{H_{0}^1(\Om)}
\newcommand{\Honet}{\Hone \times \Hone}

\newcommand{\HLap}{H(\Delta, \Om)}
\newcommand{\HzLap}{H_0(\Delta,\Om)}
\newcommand{\Hzz}{H_{00}(\Delta,\Om)}
\newcommand{\HzzLap}{\Hzz}

\newcommand{\Linfty}{L^{\infty}(\Om)}
\newcommand{\Lone}{L^1(\Om)}
\newcommand{\Loloc}{L^1_{loc}}
\newcommand{\Lp}{L^p (\Om)}
\newcommand{\Lq}{L^q (\Om)}
\newcommand{\Lt}{L^2 (\Om)}

\newcommand{\Woo}{W^{1,1}(\Om)}
\newcommand{\Wool}{W^{1,1}_{loc}(\Om)}
\newcommand{\Wop}{W^{1,p}(\Om)}
\newcommand{\Wzop}{W_0^{1,p}(\Om)}

% spaces of smooth and Harmonic functions

\newcommand{\BLap}{\Bc(\Om)}
\newcommand{\BzLap}{\Bc_0(\Om)}
\newcommand{\Cco}{C_c^1 (\Om)}
\newcommand{\Ccto}{C_c^2 (\Om)}
\newcommand{\Ccty}{C_c^{\infty} (\Om)}

\newcommand{\Harm}{\Hc (\Om)}

\newcommand{\Hharm}{{\Hc}^{1/2}(\Om)}
\newcommand{\Hsharm}{{\Hc}^{s+1/2}(\Om)}

% spaces of vector fields.
\newcommand{\CuOm}{Curl(\Om)}
\newcommand{\CuoOm}{Curl^1(\Om)}
\newcommand{\CuzOm}{Curl_0(\Om)}
\newcommand{\CuzoOm}{Curl_0^1(\Om)}
\newcommand{\CuGam}{Curl_{\Gamnu}(\Om)}

\newcommand{\GOm}{G(\Om)}
\newcommand{\GoOm}{G^1(\Om)}
\newcommand{\GpOm}{G^p(\Om)}
\newcommand{\GzOm}{G_0(\Om)}
\newcommand{\GzoOm}{G_0^1(\Om)}

\newcommand{\GradGam}{G_{\Gamtau}(\Om)}

\newcommand{\Harmf}{\Hc^0(\Om,\Rto)}
\newcommand{\HarmC}{\Hc \Cc(\Om)}
\newcommand{\HarmG}{\Hc\Gc(\Om)}
\newcommand{\Harmto}{\Hc(\Om, \Rto)}
\newcommand{\HarmtoZ}{\Hc^0(\Om, \Rto)}

\newcommand{\Hbdiv}{H_{\pal}(\Div, \Om)}
\newcommand{\Hbcurl}{H_{\pal}(\curl, \Om)}
\newcommand{\Hcurl}{H(\curl, \Om)}
\newcommand{\Hdiv}{H(\Div, \Om)}
\newcommand{\HDC}{H_{DC}(\Om)}
\newcommand{\HDCz}{H_{DC0}(\Om)}

\newcommand{\Loto}{L^1(\Om; \Rto)}
\newcommand{\Locto}{L^1_{loc}(\Om; \Rto)}
\newcommand{\Lpto}{L^p(\Om; \Rto)}
\newcommand{\Ltto}{L^2 (\Om; \Rto)}

\newcommand{\Ncurl}{N(\curl)}
\newcommand{\Ndiv}{N(\Div)}
\newcommand{\Nbcurl}{N_{\pal}(\curl)}
\newcommand{\Nbdiv}{N_{\pal}(\Div)}

% boundary or trace  spaces

\newcommand{\Hhby}{H^{1/2}(\bdy)}
\newcommand{\Hmhby}{H^{-1/2}(\bdy)}

\newcommand{\Hmsby}{H^{-s}(\bdy)}
\newcommand{\Hsby}{H^s(\bdy)}

\newcommand{\Lpby}{L^p (\bdy, ds)}
\newcommand{\Lqby}{L^q (\bdy, ds)}
\newcommand{\Ltby}{L^2 (\bdy, ds)}

\newcommand{\Wocurl}{W^1(\Curl, \Om)}
\newcommand{\Wodiv}{W^1(\Div, \Om)}
\newcommand{\WoDC}{W^1_{DC}( \Om)}

\newcommand{\Wpcurl}{W^p(\Curl, \Om)}
\newcommand{\Wpdiv}{W^p(\Div, \Om)}
\newcommand{\WpDC}{W^p_{DC}( \Om)}

% other mathematical macros
% line 300 approx

\newcommand{\bfhj}{\hbf^{(j)}}
\newcommand{\bfkj}{\kbf^{(j)}}
\newcommand{\Calph}{C^{\alpha}}

\newcommand{\Ccv}{\Cc_v}

\newcommand{\chihat}{\hat{\chi}}
\newcommand{\Cinfty}{C^{\infty}}
\newcommand{\Cinty}{C^{\infty} }

\newcommand{\delj}{\delta_j}
\newcommand{\delone}{\delta_1}
\newcommand{\delz}{\delta_0}

\newcommand{\Ecv}{\Ec_v}
\newcommand{\etanu}{\eta_{\nu}}
\newcommand{\etatau}{\eta_{\tau}}
\newcommand{\etahj}{\hat{\eta}_j}

\newcommand{\Gcq}{\Gc_q}
\newcommand{\lamgam}{\lambda_{\Gamma}}
\newcommand{\lamogam}{\lambda_1(\Gamma)}
\newcommand{\lamm}{\lambda_m}
\newcommand{\lamone}{\lambda_1}
\newcommand{\lamnu}{\lambda_{\nu} }
\newcommand{\lamtau}{\lambda_{\tau}}

\newcommand{\mut}{\tilde{\mu}}
\newcommand{\n}[1]{\left\vert#1\right\vert}
\newcommand{\nm}[1]{\left\Vert#1\right\Vert}

\newcommand{\Ocb}{\overline{\Oc}}
\newcommand{\Omb}{\overline{\Omega}}
\newcommand{\Omo}{\Om_1}
\newcommand{\Omz}{\Om_0}
\newcommand{\opal}{\oplus_{\pal}}

\newcommand{\phit}{\tilde{\vap}}
\newcommand{\psit}{\tilde{\psi}}
\newcommand{\psihm}{\psi_{hm}}
\newcommand{\psim}{\psi^{(m)}}
\newcommand{\psip}{\psi_p}
\newcommand{\psiv}{\psi_v}
\newcommand{\psivz}{\psi_{v0}}
\newcommand{\psiz}{\psi_0}

\newcommand{\ra}{\rightarrow}
\newcommand{\ROm}{R_{\Om}}
\newcommand{\ROmM}{R_{\Om M}}

\newcommand{\toxto}{2 \times 2}
\newcommand{\vapbar}{\overline{\vap}}
\newcommand{\vaph}{\vap_h}
\newcommand{\vaphm}{\vap_{hm}}
\newcommand{\vapp}{\vap_p}
\newcommand{\vapv}{\vap_v}
\newcommand{\vapvz}{\vap_{v0}}
\newcommand{\vapz}{\vap_0}

\newcommand{\vbfm}{\vbf^{(m)}}
\newcommand{\vnu}{\vbf \cdot \nu}
\newcommand{\vflux}{v_{\nu}}
\newcommand{\vt}{\tilde{v}}
\newcommand{\vtan}{v_{\tau}}
\newcommand{\vtau}{\vbf \cdot \tau}

\newcommand{\wflux}{\wbf \cdot \nu}
\newcommand{\wtan}{\wbf \cdot \tau}

%  The following are the theorem like structures used here.
%
\newtheorem{thm}{Theorem}[section]
\newtheorem{cor}[thm]{Corollary}
\newtheorem{cond}{Condition}
\newtheorem{lem}[thm]{Lemma}
\newtheorem{prop}[thm]{Proposition}
%
%

% line 420
\title[ Planar Div-curl Problems]
{Bounds  and Regularity of Solutions of  \\  Planar  Div-curl  Problems }     %line420
  
\author[Auchmuty]{Giles Auchmuty}
\address{Department of Mathematics, University of Houston, 
Houston, Tx 77204-3008 USA}
\email{auchmuty@uh.edu}

\thanks{The author gratefully acknowledges research support by NSF award DMS 11008754. 
  }

\date{May 2, 2016}

\begin{abstract} 
New 2-norm bounds are described for the least energy solutions of planar div-curl boundary value problems on bounded 
regions in space. 
Prescribed flux, tangential or mixed flux and tangential boundary conditions are treated.
A harmonic decomposition of planar fields is used to separate the solutions due to source terms from harmonic components 
that are determined by boundary data.
Some regularity results are described.
\end{abstract}

\maketitle

\section{Introduction}\label{s1}

This paper describes some properties of  solutions of  boundary value problems for div-curl systems on bounded regions 
$\Om \subset \Rto$.
This is a degenerate elliptic system of two equations in two unknowns where the existence and uniqueness conditions for solutions depend 
 on both properties of the  data and also the topology of the region and the boundary conditions.
Questions about well-posedness and uniqueness of solutions  were studied using variational methods in Alexander and Auchmuty  \cite{AA}.
Here  primary  attention is devoted to  finding energy (2-norm) bounds on solutions, their dependence on boundary conditions  and associated regularity results.

In particular some different decompositions of the  fields will be used to obtain different and better energy inequalities and regularity results.
 Prescribed flux   prescribed tangential and mixed flux and tangential  boundary conditions will be studied.
 The analysis is based on a "harmonic decomposition", which differs from the Hodge-Weyl   decompositions used in \cite{AA}.
 In this decomposition, the potentials associated with sources are solutions of zero-Dirichlet problems for Poisson's equation.
 Then the harmonic component is determined as solutions of homogeneous equations with nontrivial Neumann boundary data.

The energy (2-norm) bounds obtained here hold for the least energy solutions when the boundary value problem is underdetermined.
 To simplify the regularity statements, and also some estimates, only the simplest div-curl system is treated.
  A number of other authors have regard  this system is a prototype for degenerate elliptic systems,  while K.O. Friedrichs \cite{Fr}
has  called an inequality that bounds the the energy of a field by norms of  its divergence and curl, the {\it  main inequality   of 
vector analysis}. 
 This system has been used to  model many different situations in both fluid mechanics and electromagnetic field theories.
Many of the results  described here may be generalized to the case of general elliptic coefficient matrix $E(x)$, as used in \cite{AA}
 using standard    assumptions and techniques.
 
 The bounds described here are quite different to those in papers such as that of Krizek and Neittaanmaki \cite{KN} or
 the Schauder estimates  of Bolik and von Wahl \cite{BvW}   for example.
 Since a major interest is the dependence of solutions on boundary data, this analysis is also quite different to the work of Brezis
 and Bourgain \cite{BB} who studied these problems in $\RN$ or with periodic boundary conditions.
  A sophisticated analysis of the prescribed flux and tangential boundary   problems has been given by Mitrea \cite{Mi}.
 She used layer potentials and Besov spaces to study these problems on bounded regions with Lipschitz boundaries.
 There also has been considerable work on the numerical analysis and simulation of solutions of these problems by a large number of authors.  See Bramble and Pasciak \cite{BP} for example, or Monk \cite{Mon} for an overview.    

In section 3, some results about the regularity of orthogonal projections of $L^2-$vector fields are described. 
It is shown that, when a field is smooth on an open subset of $\Om$, so are its potentials.
(The function that is often called a stream function will be called a potential here.)
Also the  class of harmonic vector fields that can be represented by conjugate harmonic functions is characterized.

In section 4, some orthogonal decompositions of the space $\Ltto$  the standard inner product  are described.  
Readers can look at page 314 of Dautray and Lions \cite{DL3} to see the possible orthogonal representations  available.
While that  diagram is for 3-dimensional regions, essentially the  same analysis holds for planar regions.
The only simplification for 2-dimensional problems is that the dimensions of the spaces of special harmonic fields must be 
equal in 2-d while they can vary in 3-d.
Theorem 4.2 here is the crucial theorem about the existence and regularity of representations of irrotational or solenoidal 
vector fields in $\Lpto$ using  potentials.
It should be commented that the non-uniqueness of the Helmholtz decomposition of vector fields has been a source of many problems  in applications as many of the   commonly used splittings  are not described  by projections.
As a corollary a proof that harmonic vector fields must be $C^{\infty}$ on $\Om$ is proved. 
This generalizes  Weyl's lemma for harmonic functions. 

Section 5 describes some  properties of scalar Laplacian boundary value problems are collected for later use.
Some properties of the solution of zero-Dirichlet boundary value problem for Poisson's equation that are  not readily accessible
are described.
One consequence is that a harmonic decomposition that applies to $L^1$ and $\Loloc$ fields  is obtained.
A different Hilbert space $\HzLap$ is introduced so that results may be given when the boundary  $\bdy$ is not 
necessarily  $C^1$.
Explicit formulae, and estimates, for solutions of Neumann, and other, harmonic boundary value problems  in terms of Steklov eigenvalues and eigenfunctions  follow from the author's work in cite{AuH} and are needed to obtain results about 
 the dependence of solutions on boundary data.

Sections 6 and 7 describe results about the least norm solutions of the prescribed flux and prescribed tangent div-curl boundary
value problem respectively when the necessary compatibility conditions hold.
When each of the data is $L^2$,  solution estimates depending on the principal Dirichlet and Stekov eigenvalues 
$\lamone$ and $\delone$ of the Laplacian on $\Om$.  
These estimates are sharp.
When $\Om$ is not simply connected these solutions are not  unique and  the extra information required for well-posedness
were  studied in \cite{AA}. 
% In section 8 the two norm of the special harmonic fields are shown   to be given by certain relative capacities.  

When mixed tangential and normal boundary conditions are imposed on $\bdy$, no compatibility conditions on the  data
are required  for the existence of solutions. 
Under natural assumptions on the data, it is shown how the solutions may be represented using two potentials 
and estimates of these solutions are found in terms of some different eigenvalues associated with the Laplacian on $\Om$.

This paper aims to provide a self-contained description of some basic results about these problems in a manner that can be used by 
numerical analysts and others interested in the approximation, and properties,  of solutions.
Thus some of the results here are variants of results known, often in much greater generality, to researchers in linear elliptic
boundary value problems. 

  % % % % 
 % % % % %
\vspace{1.5em}

%\newpage
% ####### Section 2 ######
% This is line 530 approx
\section{Definitions and Notation.}  \label{s2}

In this paper, standard definitions as given in Evans text \cite{Ev}  or Attouch, Butazzo and  Michaille \cite{ABM} will generally be used.  - specialized to $\Rto$ since this paper  only treats planar problems. 
Cartesian coordinates  $x \eqs (x_1,x_2)$ will be used and Euclidean norms and inner products are denoted by $\, |.| \, $ and 
$ x \cdot y $. 
 A region  is a non-empty, connected, open subset of $\Rto$. Its closure  is denoted  $\Omb$ and its boundary is   
   $\bdy \deq \Omb\setminus \Om$. 
   Often the position vector $x$ is omitted in formulae for functions and fields and equality should be interpreted as holding  a.e. 
  $\wrt$ 2-dimensional  Lebesgue  measure $\dtx = dx_1 \, dx_2 $ on $\Om$.

A number of the results here  depend on the differential topology of the  region $\Om$. 
 A curve in the plane is said to be a simple Lipschitz loop if it is a closed, non-self-intersecting curve with  at least two distinct
 points and  a uniformly  Lipschitz parametrization.   
 Such loops will be compact and have finite, nonzero, length.
Arc-length will be denoted $s(.)$ and our standard assumption is
 
\noindent{\bf Condition B1. }  {\it $\Om$ is a bounded region in  $\Rto$ with boundary  $\bdy$  the union of a finite number of 
disjoint simple Lipschitz  loops $\{ \Gamj : 0 \leq j \leq J\}$.}

Here $\Gamz$ will always be the exterior loop and the other $\Gamj$ will enclose {\it holes} in the region $\Om$
The interior region to the loop $\Gamz$  defined by  the Jordan curve theorem  will be denoted $\Omz$.
When $J=0, \  \Om$ is said to be simply connected and then $\Omz = \Om$.

The outward unit normal to a region at a point on the boundary is denoted $\nu(z) = (\nu_1(z), \nu_2(z))$.
Then  $\tau(z)  := (-\nu_2(z), \nu_1(z))$ is the positively oriented unit tangent vector at a point $ z \in \bdy$. 
   $\nu, \tau $ are defined $s \, a.e.$  on $\bdy$ when (B1) holds. 

In this paper, all functions are assumed to be at least $\Loloc$  and derivatives will be taken in a weak sense. 
The spaces $\Wop, \Wzop$ are defined as usual for $ p \in [1,\infty]$ with standard  norms  denoted by $\| . \|_{1,p}$. 
When $p = 2$ the spaces will also  be denoted $\Hone, \Hzone$.

 When $\Om$ is bounded,   the trace of  Lipschitz continuous functions on  $\Omb$ restricted to  $\bdy$ is again Lipschitz continuous.
 The extension of this linear mapping  is a continuous linear mapping of $W^{1,p}(\Om)$ to $\Lpby$ for all $p \in [1,\infty]$ when (B1) holds.
 See \cite{EG}, Section 4.2 for details. 
From Morrey's theorem,  $\gamma$ maps $\Wop$ into $\Calph(\bdy) $ when $p > 2$ and $\alpha = 1-2/p$.  
di Benedetto \cite{DiB} proposition 18.1 shows that when $\vap \in \Hone$ then $\gamphi \in \Lqby$ for all $q \in [1,\infty)$.
Also if $\vap \in \Wop$ with $p \in [1,2)$ then $\gamphi \in \Lqby$ for all $q \in [1,p_T]$ with $p_T = p/(2-p)$  under stronger regularity   conditions on the boundary.  

The region $\Om$ is said to satisfy a {\it compact trace theorem} provided  the  trace mapping $\gamma : \Hone \ra \Ltby$ is compact. 
% Often $u$ will be used in place of $\gamma( u)$ for the trace of a function on $\bdy$.
Theorem 1.5.1.10 of Grisvard  \cite{Gr} proves an inequality that implies the compact trace theorem when $\bdy$  satisfies (B1).  

We will generally  use  the following  equivalent inner product on $\Hone$
\begin{equation}\label{ip2}
{[\vap , \psi]}_{\pal}  \deq \IOm \  \gradphi  \cdot \nabla \psi \  \dtx \pls  \Iby \  \gamphi \; \gampsi \ ds
\end{equation}
The associated  norm is denoted $ {\nm{\vap}}_{\pal}$.
   The proof that this  norm is equivalent to the usual $(1,2)-$norm on $\Hone$ when  (B1) holds is Corollary 6.2 of \cite{AuSE} 
  and also is part of theorem 21A of \cite{Z2A}. 
  The inner product on $\Hzone$ is the restriction of this inner product.  

When  $\Om$ satisfies (B1), then the {\it Gauss-Green} theorem  holds in the forms
\begin{eqnarray} 
\IOm D_j \vap(x) \dtx & \eqs & \Iby \, \gamma( \vap)(z) \ \nu_j(z) \, ds(z) \ \foral \vap \in \Woo \xand \\
\IOm   \vap(x) D_j \, \psi(x) \  dx & = &  \Iby \gamma(\vap) \, \gamma( \psi)  \, \nu_j  \dsg \, -  \IOm   \psi(x) D_j \, \vap(x) \  dx \quad \mbox{for each} \ j
 \label{GG} 
 \end{eqnarray}
and all $\vap,\ \psi$ in $\Wop$ with $p \geq 4/3$. 
Often the trace operator will be  implicit in boundary integrals

 When $\vap \in \Woo$ is weakly differentiable, then the  {\it gradient}   and  {\it Curl} of $\vap$ are the vector fields
\beq \label{e2.5}
  \gradphi (x) := (D_1\vap(x), D_2\vap(x))  \xand  \gradp \vap(x) := (D_2 \vap(x), - D_1\vap(x)). \eeq
Here $D_j \vap$ or $\vap_{,j}$ denotes the weak j-th derivative.

A function $\rho \in \Loloc$ is defined to be the Laplacian of $\vap$ provided one has
\[ \IOm \vap \ \Lapv \dtx \eqs \IOm \ \rho \ v \ \dtx \foral v \in \Ccto \]
A function $\vap \in \Woo$ is said to be harmonic on $\Om$ provided
\beq \label{Hfn}
\IOm \gradphi \cdot  \gradchi \  \dtx \eqs 0 \  \foral \chi  \in \Ccto .  \eeq

The subspace of all harmonic functions in $\Hone$ will be denoted $\Harm$ and it is straightforward to observe that 
$\Hone = \Hzone \oplus_{\pal} \Harm$  and that $\Harm$ is isomorphic to the trace space $\Hhby$.
Later  use will be made of the analysis in \cite{AuTr} where this is described and $\pal-$orthogonal bases of the space
 $\Harm$ are found that involve the Steklov eigenfunctions of the Laplacian on $\Om$.

\vspace{2em} 
%%%%%
%%%%%
%%%%%   Line 600
  
\section{Projections and Potentials in $\Ltto$.}   \label{Projs}

Here we will first describe the representation of planar vector fields by scalar potentials $\vap, \psi$ in the form 
\beq \label{e3.1}
\vbf(x) \eqs \gradp \psi (x) \mns \gradphi(x) \ \xon \Om. \eeq
Often $\psi$ is called a stream function and a Cartesian frame on $\Om$ is used. 
This representation is generally called a Helmholtz decomposition  and many different  choices of $\vap, \psi$ have been used by scientists and engineers  for different boundary value problems.
The choice of signs in \eqref{e3.1} is commonly used in applications and also introduces some mathematical consistency. 
In this paper,  decompositions of the form \eqref{e3.1} that are defined by projections and also have  orthogonality properties will be analyzed in
some detail.

For $p \in [1,\infty], \ \Lpto$ is the space of  planar vector fields $\vbf{\it (x) = ( v_1(x), v_2(x))} $  on $\Om$ whose components  
are $L^p-$functions on $\Om$. 
  Let $\Ltto$ is the real Hilbert space of $L^2-$vector fields on $\Om$ with the inner product 
\beq \label{e3,2}
\ang{\vbf, \wbf} \deqs \IOm \vbf \cdot \wbf \ {\it \dtx}. \eeq
Throughout this paper if fields or subspaces are said to be orthogonal, without any further adjectives, this inner product is implied.

Define $\GOm, \GzOm, \CuOm, \CuzOm$ be the subspaces of  gradients and Curls with potentials $\vap$ in $\Hone, \Hzone$  respectively. 
These subspaces will first be shown to be  closed subspaces of $\Ltto$ and then some properties of the
associated orthogonal projections are obtained.
This will be done using a variational characterization of projections  based on Riesz' projection theorem as described in
section 3 of  Auchmuty \cite{AuOD}.

First consider   the approximation  of two-dimensional vector fields by gradient fields.
This becomes a  problem of minimizing the functional $\Ecv$ defined by
\beq \label{ePP1}
\Ecv (\vap) \deqs \IOm  \ [ \, |\gradphi|^2 \pls 2 \, \vbf \cdot \gradphi \, ] \ \dtx \eeq
on $\Hone, \Hzone$ respectively. 
This functional differs from  $\| \vbf \pls \gradphi \|^2$ by $\| \vbf\|^2$ so minimizing  this is equivalent to finding the 
best approximation of $\vbf$ by   gradients in $\Ltto$.

As is standard, the space $\Hone$ is replaced by the space $\Hmone$ of all potentials with mean value $\vapbar = 0$. 
The inner product on both $\Hmone$ and $\Hzone$ is taken as $\ang{\vap, \chi}_{\nabla} \deqs \ang{\gradphi, \gradchi}$. 

The results about these variational principles may be summarized as follows.
\btm \label{TP1}
Assume $\Om$ obeys (B1) and $v \in \Ltto$. 
Then there is a unique $\vapv \in \Hmone$ that minimizes $\Ecv$ on $\Hmone$ and it  satisfies
\beq \label{e3.3}
\IOm \, (\gradphi \pls \vbf) \cdot \gradchi \ {\rm \dtx \eqs 0 \foral  \chi \in \Hone.} \eeq
Moreover (i)  \quad $\| \, \gradphi_v \| \leqs \| \vbf \|$ and  (ii) \quad  if $\vbf = - \, \gradpsi$ then there is a constant c such that 
$\vapv + \psi \equiv c$ on $\Om$.
\etm    \bpf
When (B1) holds there is a $\lamm > 0$ such that 
\beq \label{3.4}
\IOm \ |\gradphi|^2 \, \dtx  \geqs  \lamm \, \IOm \vap^2 \, \dtx \foral \vap \in \Hmone.   \eeq
Hence $\Ecv$ is continuous, strictly convex and coercive on $\Hmone$ so there is a unique minimizer of $\Ecv$.
This functional is Gateaux differentiable and  the minimization condition is \eqref{e3.3}. 
Choose $\chi = \gradphi_v$ in \eqref{e3.3} then (i) follows from Cauchy-Schwarz. 
(ii) follows as each function in $\Hone$ has a unique decomposition of the form $\vap = \vap_m + c$ where $\vap_m \in \Hmone$
and $c = \vapbar$. 
\epf

Define $\PG : \Ltto \to \Ltto$ by $\PG \, \vbf \deqs - \gradphi_v$. 
This result implies that   $\GOm$ is a closed subspace of $\Ltto $ from corollary  3.3 of cite{Au1}  and (ii) shows that $\PG$ is the  
projection of $\Ltto$ onto $\GOm$.
The extremality condition \eqref{e3.3} implies that $\vapv$ is a weak solution of the Neumann problem
\beq \label{e3.5}
- \, \Lapphi \eqs \Div \, \vbf \xon {\rm \Om \xand \Dnu \vap \eqs} - \, \vbf \cdot \nu {\rm \xon \bdy.} \eeq

The orthogonal complement of this projection is $\QG := I - \PG$ and is the 2d version of the Leray projection of fluid mechanics. 
Equation \eqref{e3.3} says that $\QG$ and $\PG$ are orthogonal projections. 

\btm \label{TP2}
Assume $\Om$ obeys (B1) and $v \in \Ltto$. 
Then there is a unique $\vapvz \in \Hzone$ that minimizes $\Ecv$ on $\Hzone$ and it  satisfies
\beq \label{e3.7}
\IOm \, (\gradphi \pls \vbf) \cdot \gradchi \ {\rm \dtx \eqs 0 \foral \chi \in \Hzone} \eeq
and  \quad $\| \, \gradphi_{v0} \| \leqs \| \vbf \|$. \etm    \bpf
When (B1) holds there is a $\lamone> 0$ such that 
\beq \label{e3.33}
\IOm \ |\gradphi|^2 \, \dtx  \geqs  \lamone \, \IOm \vap^2 \, \dtx \foral \vap \in \Hzone.   \eeq
Hence $\Ecv$ is continuous, strictly convex and coercive on $\Hzone$ so there is a unique minimizer of $\Ecv$.
This functional is Gateaux differentiable and  the minimization condition is \eqref{e3.7}. 
Choose $\chi = \gradphi_v$ in \eqref{e3.7} then the last part follows from Cauchy-Schwarz. 
\epf

The extremality condition \eqref{e3.3} says that $\vapvz \in \Hzone$ is a weak solution of the Dirichlet problem for
\beq \label{e3.8}
- \, \Lapphi \eqs \Div \, \vbf \xon {\rm \Om}   \eeq

Define $\PGz : \Ltto \to \Ltto$ by $\PGz \, \vbf \deqs - \,  \gradphi_{v0}$. 
This result implies that   $\GzOm$ is a closed subspace of $\Ltto $ from corollary  3.3 of \cite{AuOD}  and (ii) shows that $\PGz$ is the  
projection of $\Ltto$ onto $\GzOm$.
The complementary projection $\QGz : = I - \PGz$ is the projection onto the null space of the  divergence operator - see theorem 
\ref{T4.2} of the  next section.

These characterizations of these projections  allow the proof that they preserve interior regularity.
A vector field $\vbf$ is said to be $H^m$ on an open subset $\Oc$ provided each component $v_j$ is of class $H^m$ on $\Oc$. 
Specifically the following holds. 

\btm \label{TR1}
Suppose $\Om$ obeys (B1) and $\Oc$ is  open   with $\Ocb \subset  \Om$. 
If $\vbf \in {\rm H^m(\Oc)}$ with $m \geq 1$,  then so are $\PG \vbf, \, {\it \PGz} \vbf, \, {\it \QG } \vbf, $ and $\QGz \vbf$ .
\etm

\vspace{-1em}
\bpf
When $\vbf \in {\it H^m}$  on $\Oc$, then $\Div \vbf$ is $H^{m-1}$ and thus $\vapv$ is $H^{m+1}$ from standard elliptic 
regularity results for solutions of \eqref{e3.5}   as in Evans, \cite{Ev} chapter 6 or elsewhere. 
Hence the gradient is $H^m$ so the results hold for $\PG \vbf$ and $\QG \vbf$ 
The result for $\PGz \vbf, \QGz \vbf$ is proved in the same way since the potentials now are  solutions of \eqref{e3.8}.   
\epf

Analogous analyses hold for  projections onto spaces of Curls. 
First note that 
\[ \| \vbf \mns \gradp \psi \|^2 \mns \| \vbf \|^2 \eqs \IOm \,  [ \, |\gradpsi |^2 \mns 2 \, \vbf \wedge \gradpsi \, ] \ {\it \dtx} \]
where $ \wedge$ denotes the 2d vector product. 
Consider the variational problems of minimizing  the functional
\beq \label{VFl2}
\Ccv (\vap) \deqs \IOm  \ [ \, |\gradphi|^2 \mns 2 \, \vbf \wedge\gradphi \, ] \ {\it \dtx} \eeq
on $\Hone, \Hzone$ respectively. 
Results about these variational principles may be summarized as follows.

\btm \label{TP3}
Assume $\Om$ obeys (B1) and $v \in \Ltto$. 
Then there is a unique $\psiv \in \Hmone$ that minimizes $\Ccv$ on $\Hmone$ and it  satisfies
\beq \label{e3.6}
\IOm \, (\gradp \psi \mns \vbf) \cdot \gradp \chi \ {\rm \dtx \eqs 0 \foral  \chi \in \Hone.} \eeq
Then (i) \quad $\| \, \gradp \psi_v \| \leqs \| \vbf \|$,  \quad  (ii) \quad  if $\vbf = \gradp \psi$ then there is a constant c such that 
$\psiv - \psi \equiv c$ on $\Om$ and (iii) \quad
If  $\Oc$ is an open subset of $\Om$ with $\Ocb \subset \Om$ and $\vbf \in {\rm H^m(\Oc)}$, then $\gradp \psiv \in H^m(\Oc)$.   \etm 
\vspace{-1em}
\bpf
This proof just involves   appropriate modifications to those of theorems \ref{TP1}, and \ref{TR1}
\epf

The extremality condition \eqref{e3.6} says that $\psiv$ is a weak solution of the Neumann problem
\beq \label{e3.75}
- \, \Lappsi \eqs \curl\, \vbf \xon {\rm \Om \xand \Dnu \psi \eqs} - \vbf \cdot \tau {\rm \xon \bdy.} \eeq

Define $\PC : \Ltto \to \Ltto$ by $\PC \, \vbf \deqs \gradp \psi_v$. 
The theorem  implies that   $\CuOm$ is a closed subspace of $\Ltto $ from corollary  3.3 of cite{Au1}  and (b) shows that $\PC$ is the  
projection of $\Ltto$ onto $\CuOm$.
The orthogonal complement of this projection is $\QC := I - \PC$ and  \eqref{e3.3} says that $\QC$ and $\PC$ are orthogonal projections. 

Similarly, the problem of minimizing the functional $\Ccv$ on $\Hzone$ has solutions that satisfy the following. 
The proof is similar to that of theorems \ref{TP2} and \ref{TR1}. 
\btm \label{TP4}
Assume $\Om$ obeys (B1) and $v \in \Ltto$. 
Then there is a unique $\psivz \in \Hzone$ that minimizes $\Ccv$ on $\Hzone$ and it satisfies
\beq \label{e3.9}
\IOm \, (\gradp \psi \mns \vbf) \cdot \gradp \chi \ {\rm \dtx \eqs 0 \foral  \chi \in \Hzone.}  \eeq
Thus (i) \quad $\| \, \gradp \psi_{v0} \| \leqs \| \vbf \|, $ (ii) \quad  if $\vbf = \gradp \psi$  with $\psi \in \Hzone$ then  
$\psivz = \psi$  and (iii) \quad
if  $\Oc$ is an open subset of $\Om$ with $\Ocb \subset \Om$ and $\vbf \in {\rm H^m(\Oc)}$, then $\gradp \psivz \in H^m(\Oc)$.
\etm 

Define $\PCz : \Ltto \to \Ltto$ by $\PCz \, \vbf \deqs \gradp \psi_{v0}$. 
This result implies that   $\CuzOm$ is a closed subspace of $\Ltto $ from corollary  3.3 of \cite{AuOD}  and (iii) shows that $\PCz$ is the  
projection of $\Ltto$ onto $\CuzOm$.
The extremality condition \eqref{e3.9} says that $\psivz \in \Hzone$ is a weak solution of the Dirichlet problem for
\beq \label{e3.10}
- \, \Lapphi \eqs \curl\, \vbf \xon {\rm \Om} .  \eeq

The complementary projection $\QCz : = I - \PCz$ is the projection onto the null space of the curl operator  as described in  the next section.
These results may be combined to yield the following result that has  been central in the study of 2-dimensional perfect fluids
and much classical study of vector fields.

\btm \label{THarm}
 Assume $\Om$ satisfies (B1) and $\vbf \in {\it \GOm \cap \CuOm}$ then the potentials $\vapv, \psiv$  are conjugate 
 harmonic functions on $\Om$ and $\vbf$ is $\Cinfty$ on $\Om$. \etm
 \bpf
 The assumption is that there are functions $\vap_v, \psi_v$ such that  $\vbf \eqs  {\it \gradphi_v \eqs \gradp \psi_v}$ on $\Om$.  
 These are the Cauchy-Riemann equations. 
 Then \eqref{Hfn} follows for each of $\vapv, \psiv$ upon using Gauss Green and the commutativity of weak derivatives. 
 The potentials are harmonic functions so they are $\Cinfty$ and thus $\vbf$ is.  
  \epf

% Should introduce regularity results and reductions to harmonic vector fields here?
% WHat happens when $p = 1$?

%%%%%
%%%%%
%%%%%
%\newpage

  % ###   Section 4 ###
  % This is line 920 approx
\vspace{2em}
  \section{Div, Curl and Orthogonality of Planar Vector Fields} \label{s4} 
 
A basic question for these fields is how the projections onto the spaces of gradients and Curls defined above are related to
the vectorial operators $\Div$ and $ \curl$ on Sobolev-type spaces  of vector fields?
In this section such spaces are described and some orthogonality results obtained.
 
 The  curl of  a vector field  $\vbf \in {\rm \Locto}$ is a function $\omega \in \Loloc(\Om)$, (or possibly 
  a distribution) that satisfies 
\beq \label{e4.1}
 \IOm \gradp \psi \cdot \vbf \  {\it \dtx} \eqs \IOm \omega\, \psi \ {\it \dtx}  \foral {\it \psi \in \Cco }. \eeq
In this case we write $\curl \vbf \deqs \omega$.
 Similarly the  divergence of  $\vbf  {\rm \in  \Locto}$   is defined to be the function 
 $\rho \in \Loloc(\Om)$ provided $\rho$ satisfies 
\beq \label{e4.2}
 \IOm \gradphi \cdot \vbf\ {\it \dtx} \eqs  - \IOm \rho \, \vap \ {\it \dtx}  \foral {\it \vap \in \Cco. }\eeq
When this holds we write    $\Div \vbf \deqs \rho$.  

In the following only fields whose curl and div are locally  integrable functions  on $\Om$ will be studied.
Note these definitions do not require that the individual components of the derivative matrix are finite, or even defined.
When the components of a planar vector field $\vbf$ are in $\Woo$, then the derivative of the field is the matrix valued 
function $D \,\vbf(x) := (v_{j,k})$ whose entries are $L^1$ functions on $\Om$. 
It is straightforward to verify that then
\beq \label{e4.3}
\Div \vbf \eqs v_{1,1} + v_{2,2} \quad  \xand \quad \curl \vbf \eqs v_{2,1} - v_{1,2}. \eeq

A field $\vbf \in {\rm \Locto}$ is {\it irrotational}, or {\it solenoidal}, respectively provided
\beq \label{e4.5}
  \IOm \gradp \phi \cdot \vbf \  {\it \dtx} \eqs 0  \xor  \IOm \gradphi \cdot \vbf\ {\it \dtx} \eqs  0 \foral {\it \vap \in \Cco }. \eeq
 A field $\vbf \in {\rm \Locto}$ is {\it harmonic}  if it is both irrotational and solenoidal on $\Om$.
Let $\Harmf$ be the closed subspace of all harmonic vector fields in $\Ltto$. 
Observe that the space $H_{GC}(\Om) := \GOm \cap \CuOm$ of  vector fields that are both gradients and curls of $H^1-$ functions 
is a space of harmonic fields. 

 Let $\Ncurl, \Ndiv$ be the subspaces of irrotational, solenoidal vector fields in $\Ltto$ respectively.  
Note  that fields in $\GOm $ are in $\Ncurl$ and fields in $\CuOm$ are in $\Ndiv$ from the commutativity of weak differentiation.
  A first orthogonal decomposition result is the following CGH decomposition - which is independent of the differential topology of 
  $\Om$.
  
\btm \label{TDec1}
 Assume $\Om$ satisfies (B1), then 
\vspace{-0.5em}
\begin{eqnarray}   \label{e4.6}
   (a) \qquad \  \Ltto & \eqs & \ \CuzOm \opluss \Ncurl  \eqs \GzOm \opluss \Ndiv,  \xand   \\  
  (b) \qquad \ \Ltto & \eqs &  \ \CuzOm \opluss  \GzOm \opluss  \Harmf.        \label{CGH}
  \end{eqnarray}       \etm
  
  \vspace{-1.2em}
\bpf
(a) follows from the definition of $\Ncurl$ and $\Ndiv$ since $\Cco$ is dense in $\Hzone$. 
Then (b) follows as $\vbf \in {\rm \Harmf}$ iff  it is orthogonal to both $\CuzOm$ and $\GzOm$.
\epf

This theorem implies that  the projections $ \QGz, \QCz$ of the preceding section are the projections onto the 
subspaces $\Ndiv, \Ncurl$ respectively  since $\PGz, \PCz$ are the projections onto their orthogonal complements.

The splitting of \eqref{CGH} will be called the {\it  harmonic decomposition}  and will be the primary representation 
used from now on in this paper. 
It is different to the usual Hodge-Weyl decompositions where the zero boundary conditions are imposed on only one of the
potentials.  
 The following  common physical sign convention will be used. 
\beq \label{CGHR}
\vbf \eqs \gradp \psi \mns \gradphi \pls \hbf   {\rm  \xwith \psi,  \ \vap  \in \Hzone \xand } \hbf \in {\rm \Harmf.}  \eeq

Let $\Oc := I_1 \times I_2$ be an open rectangle in $\Rto$.
$\Pslemma$ provides  explicit formulae for the potentials $\vapp, \psip$ of irrotational and solenoidal $C^1-$fields on $\Om$.
Namely given a point $P \in \Oc$ and a piecewise $C^1-$ curve $\Gamx$ joining $P$ to $x = (x_1,x_2) \in \Oc$, then
\beq \label{e4.8}
\vapp(x) \deqs \int_{\Gamx} \ v_1 \, dx_1 + v_2 \, dx_2  \xand \psip(x) \deqs \int_{\Gamx} \ v_1 \, dx_2 -  v_2 \, dx_1 \eeq
are well-defined $C^1-$ functions on $\Om$.
When $\vbf$ is irrotational, $\gradphi_p = \vbf$ and when $\vbf$ is solenoidal then $\curlpsi_p = \vbf$ on $\Oc$. 
 See Dautray and Lions, \cite{DL3},  Chapter IX, section 1, lemma 3 for a proof in   the case where $\Om$ is a block in $\Rte$. 
 The proof there is easily modified for this 2 dimensional case.

The line integrals in \eqref{e4.8} are not well-defined when the field $\vbf$ is only $L^p$ on $\Om$.
Nevertheless, potentials may be proved to exist using a density argument.
The following result  is known for such fields with $p = 2$ and $\Oc \subset \Rte$; see Girault - Raviart \cite{GR} or 
Monk, \cite{Mon} theorem 3.37.

\btm \label{T4.2}
 Assume $\Om$ satisfies (B1) and $p \in [1,\infty)$.
 If $\vbf \, {\it \in \Lpto}$ is irrotational then there is a $\vap \in \Wop$ such that $\vbf = \gradphi $ on $\Om$.
  If $\vbf \, {\it \in \Lpto}$ is solenoidal then there is a $\psi \in \Wop$ such that $\vbf = \gradp \psi $ on $\Om$.  \etm
\bpf 
First assume that $\Om$ is convex, then Poincar\'{e}'s lemma implies  this result holds when $\vbf$ is $C^1$ on $\Om$.
To prove this holds for any $L^2$ field introduce a $C^1-$ mollifier $\Phi$ and consider fields on the open convex neighborhood 
$\Omo$ of points within distance 1 of $\Om$. 
The sequence of $C^1-$ fields $\vbfm$ defined by convolution $\vbfm :=  {\rm \Phi_m \star}  \vbf$   converges to the zero 
extension of $\vbf$  to $\Omo$ in $\Lt$. 
Each of these $\vbfm = \gradp \psim $ on $\Omo$ from Poincar\'{e}'s lemma. 
Normalize the $\psim$ to have mean value zero.
Since these fields are a Cauchy sequence in $\Lt$, the $\psim$ are Cauchy in $H^1_m(\Omo) $, so they converge to a limit $\psit$.
Taking limits, $\gradp \psit$ is the zero extension of $\vbf$ to $\Omo$, and the result follows by considering the restriction 
to $\Om$.
A similar proof works for the second part, using the second part of the classical Poincar\'{e} lemma.

When $\Om$ is not convex, then choose $\Omo$ in the above proof be the neighborhood of distance 1 from the convex hull of $\Om$.
Then the same arguments yield the statement of the theorem.
\epf

Note that the preceding proof extends to 3-dimensional vector fields and regions, with the usual modifications, as the construction 
 of Poincar\'{e}'s lemma is  valid there - and the other ingredients are independent of dimension. 
A corollary is the following vector-valued version of Weyl's lemma, that also extends to 3-d vector fields.
\begin{cor} \label{C4.3}
\quad  Assume $\Om$ satisfies (B1) and $p \in [1,\infty)$.
 If $\vbf \, {\it \in \Lpto}$ is a harmonic vector field then it is $C^{\infty}$ on $\Om$. 
 \end{cor}
 \bpf
Since $\vbf$ is irrotational, there is a $\vap \in \Wop$ such that $\vbf = \gradphi$ on $\Om$.
As $\vbf$ also is  solenoidal,  $\vap$ is a weak solution of Laplace's equation. 
Thus, from Weyl's lemma, $\vap$ is $C^{\infty}$ on $\Om$, and thus  $\vbf$ is also.   \epf

It should be noted that the above results  do not require any topological conditions on the region $\Om$.
The theorem implies   that $\GOm^{\perp} \subset \CuOm$ and $\CuOm^{\perp} \subset \GOm$ for any region $\Om$ satisfying (B1)
- and it is well-known that these are strict inclusions when $\Om$ is not simply connected.

\vspace{1em}

%%%%%
%%%%%   line 920
%%%%%

% The following is section 5
\section{Div-Curl and Laplacian  Boundary Value Problems}   \label{s5}

The  div-curl boundary value problem is to find a vector field $\vbf$ defined on a bounded region $\Om \subset \Rto$
that satisfies
\beq   \label{eDC.1}
\Div \, \vbf(x)  \eqs \rho (x) \qquad \mbox{and} \quad  \curl \, \vbf (x) \eqs \omega {\it (x) \qquad \mbox{for}\quad x \in \Om }
\eeq
subject to prescribed boundary conditions on $\bdy$. 

Generally either the normal component $\vbf \cdot \nu$, or the tangential component $\vbf \cdot \tau$,  of the field at  the  boundary
are prescribed in applications. 
When the normal component is prescribed everywhere on the boundary we have  a {\it normal Div-Curl boundary value problem}
 that will be analyzed in the next section. 
Problems where the tangential component is prescribed everywhere are called the {\it tangential Div-Curl boundary value problem}
 and are studied in section \ref{Tang}.
 When normal components are prescribed on part of the boundary and tangential components on the complementary subset, it 
 will be called a {\it mixed Div-Curl boundary value problem}. 

To obtain bounds and regularity results for these problems,  some results about Laplacian boundary value problems on regions obeying (B1) are required.
Although this is a standard example of a second order elliptic boundary problem, the author has not found many of
these results in the literature - so they are proved here for completeness.
Stronger regularity results are well-known when the boundary $\bdy$ is $C^k$ with $k \geq 1$ or solutions are sought 
 in various Schauder spaces. 
For many physical and numerical problems, however, it is desirable to have results on Lipschitz regions that allow
 for  "corners."

The solutions of these boundary value problems will be found by introducing appropriate potentials $\vapz, \psiz$ so that 
(b) of theorem \ref{TDec1} can be used. 
They will be solutions  of  Poisson's equation with zero Dirichlet boundary data and are characterized by variational principles.

Given $\rho \in  \Lp$, consider the problem of minimizing  the functional $\Dc$ defined by
\beq \label{e5.2}
\Dc (\vap) \deqs \IOm  \ [ \, |\gradphi|^2 \mns 2 \, \rho \, \vap \, ] \ \dtx \eeq
on $ \Hzone$. 
The essential results about this classical problem may be summarized as follows
\btm \label{T5.1}
Suppose (B1) holds, $p \in (1,\infty]$ and $\rho \in \Lp$.
Then there is a unique minimizer  $\vapz := \GD \rho$ of $\Dc$ on $\Hzone$ and it satisfies 
\beq \label{e5.1}
\IOm [ \, \gradphi \cdot \gradchi \mns \rho \chi \, ] \eqs 0 \foral \chi \in \Hzone. \eeq
The linear operator   $\GD$ is a 1-1 and  compact map of $\Lp$ into $\Hzone$. \etm
\bpf
When $\Om$ is bounded then the imbedding $i : \Hzone \to \Lq$ is continuous for all $q \in [1,\infty)$ from the Sobolev
imbedding theorem.
Thus the linear term in $\Dc$ is weakly continuous.
The existence of a unique  minimizer then holds as $\lamone > 0$  in \eqref{e3.33}, so $\Dc$ is continuous, convex and coercive
on $\Hzone$. 

$\Dc$ is G-differentiable on $\Hzone$ and the extremaility condition for a minimizer is \eqref{e5.1}
When the minimizer is denoted $\GD \rho$ it follows that $\GD$ is a linear mapping that  satisfies
\beq \label{e5.29}
\|  \gradphi \|_2 \leqs C_{p'} \ \| \rho \|_p \eeq
with $p'$ conjugate to p. Here $C_q$ is the imbedding constant for $\Hzone$ into $\Lq$.
Thus $\GD$ is continuous. It is 1-1 from the maximum principle for  harmonic functions.

 To prove $\GD$ is compact let $\{\rho_m : m \geq 1 \}$ be a weakly convergent sequence in $\Lp$ with $p > 1$.
The imbedding of $\Lp$ into $\Hnone$ is compact for $p \in (1,\infty)$ by duality to the Kondratchev theorem, 
so  the sequence $\{\rho_m \}$ is strongly convergent in $\Hnone$. 
 A standard result is that  $\GD$ is a continuous linear map of $\Hnone$ to $\Hzone$ so it is a compact linear mapping of $\Lp$ to
 $ \Hzone$ when $p > 1$ by composition.
\epf

It is worth noting that this result enables proofs of many of the results about the approximation of solutions of \eqref{e5.1}
by eigenfunction expansions in terms of the eigenfunctions of the zero-Dirichlet Laplacian eigenproblem.
It is well-known that orthonormal bases of $\Hzone$ of such eigenfunctions can be constructed.
Since $\GD$ is compact, finite rank approximations using these eigenfunctions will converge to the solution $\GD \rho$
for all $\rho$ in these $\Lp$ and this solution has the standard spectral representation arising from the spectral theorem
for compact self-adjoint maps on $\Lt.$

This result enables a generalization of the harmonic decomposition of $L^2$ fields  to fields in $\Loto$ with divergence and curl in
$L^p$ as  follows.

\btm \label{T5.2}
{\rm (Harmonic Decomposition)} \ Suppose (B1) holds, $\vbf \ {\it \in \Loloc(\Om; \Rto)}$, (or $\Loto$),  with
 $\Div \vbf, \  \curl \vbf \ {\rm \in \Lp}$   for some   $p > 1$.
Then there are $\vapz, \psiz \in \Hzone$ and a harmonic field $\hbf \ {\it \in \, \Loloc(\Om; \Rto)}$, (or $\Loto$) such that 
\beq \label{e5.23}
\vbf  \eqs \gradp \psi \mns \gradphi \pls \hbf   {\rm  \xon \Om.}  \eeq \etm
\bpf
Let  $\rho := \Div \vbf , \ \omega := \curl \vbf$ and $\vapz, \psiz$ are the associated solutions of \eqref{e5.1}.
They exist and are in $\Hzone$ from the theorem. 
Then $\hbf := \vbf -  \gradp \psi + \gradphi$ is a harmonic field that is in $\Loloc(\Om; \Rto)$, ( or $\Loto$),
 respectively when $\vbf$ is. \epf

Note that the three components in this decomposition are $L^2-$ orthogonal.
Also potentials $\vapz, \psiz$ and the harmonic field $\hbf$ here will have better regularity when stronger conditions
are imposed on the field, its divergence or curl using standard results from regularity theory.

For div-curl problems we seek potentials in the  subspace  $\HzLap$ be the subspace of $\Hzone$ of all functions 
whose Laplacians are in $\Lt$.
This is a real Hilbert space $\wrt$ the inner product
\beq \label{e5.3}
\ang{\vap, \chi}_{\Delta}  \deqs \IOm \, [ \, \Lapphi \, \Lapchi \pls \gradphi \cdot \gradchi \, ] \ \dtx    \eeq

When $\Om$ satisfies (B1) and $\bdy$ is $C^1$, then it is well known that $\HzLap = \Hzone \cap H^2(\Om)$.
See Evans chapter 8 for example.   %Check this.
Under weaker conditions on $\bdy$ such as our (B1), this need not hold - such issues have been studied by Grisvard \cite{Gr},
 Jerison, Koenig  and others.

Obviously the problem \eqref{e5.1}  has a unique solution $\vap = \GD \rho \in \HzLap$ if and only if $\rho \in \Lt$.
In this case  the following holds. 
\btm \label{T5.3}
Assume that (B1) holds and $\lamone$ is the constant in \eqref{e3.33}. 
Then the operator is $\GD$ is a homeomorphism of $\Lt$ and $\HzLap$ with $\vapz = \GD \rho$   satisfying
\beq \label{e5.8}
\|\vapz\|_2 \leqs \frac{1}{\lamone} \, \|\rho\|_2, \quad \|\gradphi_0\|_2 \leqs \frac{1}{\sqrt{\lamone}} \, \|\rho\|_2 \xand
 \|\Dnu \vapz\|_{2,\bdy}  \leqs C_0  \, \|\rho\|_2. \eeq 
 where $C_0 >  0$  depends only on $\Om$. \etm
\bpf
The first two inequalities here follow from the spectral representation of $\GD$ in terms of the Dirichlet eigenfunctions
of the Laplacian on $\Om$. 
The inequality for $\Dnu \vapz$ is theorem 3.2 in \cite{AuPK}.
\epf

Note also that if $\Oc$ is an open subset of $\Om$ with $\Ocb \subset \Om$ and $\rho$ is $H^m$ on $\Oc$, then 
the solution  $\vapz$ of \eqref{e5.1} will be of class $H^{m+1}$ on $\Oc$ from the usual interior regularity analysis. 
When $\rho, \omega \in \Lp$,  let $\vapz = \GD \rho, \ \psiz := \GD \omega$  be solutions of \eqref{e5.1} and consider
\beq \label{eDC.3}
\hbf(x) \deqs \vbf \mns \gradp \psiz \, {\rm (x) \pls \gradphiz \,(x)}. \eeq

Substituting in \eqref{eDC.1}, one sees that $\hbf$ will be a harmonic field with 
\beq \label{eDC.5}
\hbf \cdot \nu \eqs \vbf \cdot \nu  \pls {\rm \Dnu \vapz \xand} \hbf \cdot \tau \eqs \vbf \cdot \tau  \pls{\rm  \Dnu \psiz \xon \bdy.} \eeq
That is, the solvability of this div-curl system is decomposed into  zero-Dirichlet boundary value problems
involving the source terms and  a boundary value problem for a harmonic field. 
So  the following sections will concentrate on issues regarding different types of boundary value problems for harmonic 
vector fields.

Some related results about the Neumann problem for the Laplacian will also be needed later.
The harmonic components of solutions of our problems involve potentials    $\chi \in \Hone$ that satisfy
\beq \label{e5.9}
\IOm  \gradchi \cdot \gradxi \dtx  \mns \Iby  \eta \, \xi \, \dsg  \eqs 0 \foral \xi \in \Hone. \eeq
This is the weak form of Laplace's equation subject to $\Dnu \chi = \eta$ on $\bdy$.

To study the existence of solutions of this problem, consider the problem of minimizing the functional 
$\Nc: \Hone \to \R$ defined by
\beq \label{e5.11}
\Nc (\chi) \deqs \IOm  \  |\gradchi|^2 \ \dtx  \mns 2 \, \Iby  \eta \ \chi  \ \dsg \eeq
with $\eta \in \Ltby$. 
 A necessary  condition for the existence of a solution of this problem is that 
\beq \label{e5.12}
\Iby \ \eta \ \dsg \eqs 0 . \eeq

To prove the existence of a solution of this problem we need the fact that there is a $\delone > 0$ such that 
\beq \label{e5.14}
\IOm \ |\gradphi|^2 \, \dtx  \geqs  \delone \, \Iby \vap^2 \, \dsg \foral \vap \in \Hone \ \mbox{that satisfy \eqref{e5.12}}. 
\eeq
This $\delone$ is the first nonzero harmonic Steklov eigenvalue for the region $\Om$.

When (B1) and \eqref{e5.12} hold then a standard variational argument says that there is a unique minimizer 
$\chi = \Bc \, \eta$ of $\Nc$ on $\Hmone$ that satisfies \eqref{e5.9}. 

This solution operator $\Bc$ may be regarded as an integral operator that maps functions from $\Ltby$ to $\Harm \subset \Hone$.
 In particular it has a nice expression in terms of the harmonic Steklov eigenfunctions of $\Om$.
See Auchmuty \cite{AuSE} and \cite{AuH} for a discussion of the Steklov eigenproblem for the Laplacian on bounded regions
and \cite{AuPK} for further results about $\HzLap$.

A function $s_j \in \Harm$ is a Steklov eigenfunction for the Laplacian on $\Om$ provided it is a nontrivial solution of the system
\beq \label{e5.15}
\IOm \nabla s_j \cdot \nabla \xi \dtx \eqs \delta \Iby s_j \, \xi \dsg \foral \xi \in \Hone. \eeq
Here $\delta \in \R$ is the associated Steklov eigenvalue.
Let $\Lambda := \{ \delj ; j\geq 0 \}$  be the set of Steklov eigenvalues repeated according to multiplicity and with $\delta_j$ 
an increasing sequence.
The first eigenvalue is $\delz = 0$ and the corresponding eigenfunctions are constants on $\Omb$. 
It is a simple eigenvalue and the next eigenvalue is $\delone > 0$ of \eqref{e5.14}.
Normalize an associated set of Steklov eigenfunctions $\Sc := \{ s_j: j \geq 0 \}$ to be $L^2-$orthonormal on $\bdy$. Then
\beq \label{e5.16}
\IOm \nabla s_j \cdot \nabla s_k \dtx \eqs \delta_j \xwhen j = k \xand 0 \quad   \mbox{otherwise}. \eeq

Theorem 4.1 of \cite{AuTr} says that this sequence can be chosen to be an orthonormal basis of $\Ltby$ with the usual inner product.
If $\eta \in \Ltby$ satisfies the compatibility condition \eqref{e5.12} then it has the representation
\beq \label{e5.18}
\eta(z) \eqs \sum_{j=1}^{\infty} \ \etahj \, s_j(z) \xon \bdy \xwith \etahj = \Iby \eta \, s_j \dsg. \eeq

For $M \geq 1$, consider the boundary integral operators  $\Bc_M : \Ltby \to \Harm$ defined by
\beq \label{e5.19}
\Bc_M \, \eta (x) \deqs \int_{\bdy} \, B_M(x,z) \, \eta(z) \, ds \xwith B_M(x,y) \deqs \sum_{j=1}^M \, \delj^{-1} \, s_j(x) \, s_j(z) . \eeq
These are finite rank oeprators and the following SVD type representation theorem holds for the operator $\Bc$.

\btm \label{T5.4}
Assume (B1) holds, $\Lambda$ is the set of harmonic Steklov eigenvalues on $\Om$ repeated according to multiplicity
and $\Sc$  is a  $\pal-$orthogonal set of harmonic Steklov eigenfunctions and a orthonormal basis of $\Ltby$.
When $\eta \in \Ltby$,  the unique solution $ \Bc \eta$  of \eqref{e5.9} in $\Hmone$  is 
\beq \label{e5.20}
\chi(x) \eqs \Bc \eta(x) \eqs \lim_{M \to \infty}  \ \Bc_M \, \eta(x).   \eeq
 $\Bc$ is a continuous linear transformation of $\Ltby$ to $\Harm$ with  $\| \nabla \chi \|_2 \leqs  {\delone}^{-1}  \  \| \eta \|_{2, \bdy} $
 and $\chi$ is $C^{\infty}$ on $\Om$.
\etm
\bpf  
The first part of this theorem is proved in \cite{AuTr} where it is shown that $\Sc$ is an orthonormal basis of $\Ltby$.
Thus \eqref{e5.18} holds.
When $\chi$ is a solution of equation  \eqref{e5.9}, then $\chi(x) \eqs \sum_{j=1}^{\infty} \chihat_j \, s_j(x)$ on $\Om$
as $\Sc$ is a maximal orthogonal set in $\Harm$. 
Take $\xi = s_j$ in  \eqref{e5.9} then the coefficients $\chihat = \etahj / \delj$ for $j \geq 1$, and  \eqref{e5.20} holds in the 
$\pal-$norm of $\Harm$.
The function $\chi$ is $C^{\infty}$ as it is harmonic and the bound on $\| \nabla \chi \|_2$ follows from the orthogonality of $\Sc$.
 \epf

It is worth noting that the solution $\chi$ of this problem is $C^{\infty} $ on $\Om$.
It will be $H^1$ when the Neumann data is in $H^{-1/2}(\bdy)$ and more generally will be in the space $\Hc^s(\Om)$
defined as in \cite{AuH} when $\eta \in H^{s-3/2}(\bdy)$. 

%In the following sections, one of the topics of interest is what boundary conditions, together with conditions on $\Div \vbf, \curl \vbf$ imply 
%that $\vbf \in \Ltto$. 

\vspace{2em}
%%%%%
%%%%%   line 1100
%%%%%

\section{The normal Div-curl Boundary Value Problem}   \label{Norm}

The normal div-curl boundary value problem is to find a field $\vbf \, {\rm \in \Ltto}$  that solves \eqref{eDC.1} subject to
\beq   \label{BdN}
\vbf{\it (z) \cdot \nu(z)}   \eqs  {\rm  \etanu(z)    \quad   \xon \bdy.}  \eeq
with $\etanu \in \Ltby$.
From the divergence theorem,  a necessary condition for \eqref{eDC.1} - \eqref{BdN}  to have a solution is the 
compatibility condition
\beq \label{e6.2}
\int_{\Om} \rho \ \dtx \eqs \int_{\bdy}  \etanu (z) \, \dsg(z). \eeq

When the solution has the form \eqref{CGHR}, then the potentials  are solutions of \eqref{e5.1} with
$\psiz = \GD \omega, \ \vapz = \GD \rho$.
Since $\psiz \equiv 0$ on $\bdy$,  the harmonic component  satisfies
\beq \label{e6.3}
\hbf{\it (z) \cdot \nu(z) \eqs \etanu(z) \pls \Dnu \vapz (z)\quad   \xfor z \in \bdy.}  \eeq

Consider the problem of finding a gradient field that solves this problem.
If $h \eqs \gradchi$ is a solution of this problem, then $\chi$ is a harmonic function that satisfies \eqref{e5.9}
with $\eta(z) $ given by the right hand side of \eqref{e6.3}.  
From theorem \ref{T5.3}, this problem has a solution given by \eqref{e5.20} and the following result holds.

\btm \label{T6.1}
Assume (B1) holds, $\rho, \omega \in \Lt, \ \etanu \in \Ltby$ and \eqref{e6.2} holds.
Let $\vapz = \GD \rho, \ \psiz = \GD \omega$.
Then there is a unique $\chi \in \Hmone$ such that $\hbf {\it (x) \eqs \gradchi(x)}$ is a harmonic field satisfying 
\eqref{e5.9} with $\eta$  from \eqref{e6.3}. 
The field $\vbf = \gradp \psiz \mns \gradphi_0 \pls \gradchi$ is a solution of \eqref{eDC.1} - \eqref{BdN} with
\beq \label{e6.6}
\| \vbf \|_2 \leqs {\it \frac{1}{\sqrt{\lamone}} \, \left[ \ \|\rho\|_2 \pls \| \omega \|_2  \ \right] \pls \frac{1}{\sqrt{\delone}} \, 
\left[ \ \| \etanu \|_{2,\bdy}  \pls C_0 \| \rho \|_2 \ \right]. } \eeq
\etm 
\bpf
Given $\rho, \omega \in \Lt$,  theorem \ref{T5.3} yields the first two terms in the inequality \eqref{e6.6}.
Note that \eqref{e6.2} implies the compatibility condition \eqref{e5.12}, so there is a unique $\chi \in \Hmone$ that is 
harmonic on $\Om$ and satisfies the boundary condition \eqref{e6.3} from Theorem \ref{T5.4}.
The three fields in this representation of $\vbf$ are $L^2-$orthogonal so it only remains to bound $\| \gradchi \|_2$.
This bound  now follows from the last parts of theorems \ref{T5.3} and \ref{T5.4}.
\epf

It is worth noting that the constants in this inequality are best possible in that there are choices of $\rho, \omega$ and
$\etanu$ for which the right hand side equals the 2-norm of a solution of the problem. 
If $\rho, \omega \in \Lp$ for some $p > 1$, then \eqref{e5.2} implies that 
\[ \| \vbf \mns \gradchi \|_2 \leqs {\rm  C_{p'} \, [ \, \|\rho \|_p \pls \|\omega\|_p} \, ] \]

Also the regularity of the potentials $\vapz, \psiz$ here depends on the regularity of $\rho, \omega$ and the boundary $\bdy$.
They are independent of the boundary data. 
The boundary data $\etanu$  only influences the harmonic component $\gradchi$.
Moreover $\chi$ is very smooth ($C^{\infty} \cap H^1$) on $\Om$ as it is a finite energy solution of  Laplace's equation.

Quite often in the analysis of vector fields one is interested in spaces of fields such as $\Hcurl, \Hbcurl, \Hdiv  $ and $\Hbdiv$. 
These are the spaces of vector fields in $\Ltto$ that also have, respectively, $\curl \vbf, \curl \vbf$ and $\vbf \cdot \tau$, $\Div \vbf,
\Div \vbf$ and $\vbf \cdot \nu$ in $L^2$.
The above theorem helps in that the following result shows that bounds on the divergence and curl of a vector field and also of 
the boundary flux implies the field has finite energy ($L^2$ norm).

\begin{cor} \label{C6.3}
\quad Suppose (B1) holds, $\curl \, \vbf, \Div \, \vbf \, {\rm \in \Lt, \ \vbf \cdot \nu \in \Ltby}$, then $\vbf \in {\rm \Hbdiv \cap \Hcurl}$ 
and there is a  $C > 0$ such that
\beq \label{e6.8}
\| \vbf \|_2^2 \leqs {\rm C \, \left[ \ \|\, \curl \, \vbf \|_2^2 \pls \| \Div \, \vbf \,  \|_2^2   \pls   \| \vbf \cdot \nu \|_{2,\bdy}^2  \ \right]. } \eeq 
\end{cor}
\bpf
The inequality \eqref{e6.6} implies \eqref{e6.8} for an appropriate choice of C. 
Since \eqref{e6.8} holds, $\vbf$ is in both $\Hbdiv$ and $\Hcurl$ provided the region satisfies (B1).
\epf

Thus the solutions of this problem may be written as series expansions involving the Dirichlet and Steklov eigenfunctions of the Laplacian as $\GD$  and $\Bc$ from \eqref{e5.16} have representations $\wrt$ these bases of $\Hzone$ and $\Harm$ respectively. 

When $\Om$ is not  simply connected this boundary value problem has further solutions.
These was studied in \cite{AA} where the well-posed problem was shown to require the circulations around each handle  be further
specified for uniqueness. 
The above solution is the  least energy (2-norm) solution of the problem.
% Bounds on these special solutions will be described in section \ref{MultC}. 

\vspace{2em}
%%%%%
%%%%%   line 1100
%%%%%

\section{The tangential Div-curl Boundary Value Problem}   \label{Tang}

The tangential div-curl boundary value problem is to find a vector field $\vbf \,  {\rm \in \Ltto}$ that satisfies
\eqref{eDC.1} subject to
\beq   \label{BdT}
\vbf{\it (z) \cdot \tau(z)}   \eqs   \etatau(z)    \quad   \xon \bdy.  \eeq
with $\etatau\in \Ltby$.
A necessary condition, from the divergence theorem, for this problem to have a solution is that
\beq \label{e7.1}
\int_{\Om} \omega \ \dtx \eqs \int_{\bdy} \etatau (z) \, \dsg. \eeq

When the solution has the form \eqref{CGHR} , then the potentials $\vapz, \psiz$ are solutions of \eqref{e5.1}   given by 
 $\psiz = \GD \omega, \ \vapz = \GD \rho$.
Then the fact that $\vap \equiv 0$ on $\bdy$ implies that the harmonic component  satisfies
\beq \label{e7.3}
- \, \hbf{\it (z) \cdot \tau(z) \eqs  \etatau(z) \pls \Dnu \psiz (z) \quad   \xon \bdy. } \eeq

Suppose that this harmonic field is given by  $\hbf \eqs - \, \gradpchi$.
 Then $\chi$ is a harmonic function that satisfies \eqref{e5.9}  with $\eta(z) $ given by the right hand side of \eqref{e7.3}.  
From theorem \ref{T5.2}, this problem has a solution of the form \eqref{e5.20} and the following result holds.

\btm \label{T7.1}
Assume (B1) holds, $\rho, \omega \in \Lt, \ \etatau \in \Ltby$ and \eqref{e7.1} holds.
Let $\vapz = \GD \rho, \psiz = \GD \omega$.
Then there is a unique $\chi \in \Hmone$ such that $\hbf(x) \eqs - \gradpchi(x)$ is a harmonic field satisfying 
\eqref{e7.3} on $\bdy$. 
The field $\vbf = \gradp \psiz \mns \gradphiz \mns \gradpchi$ is a solution of \eqref{eDC.1} - \eqref{BdT} with
\beq \label{e7.6}
\| \vbf \|_2 \leqs \frac{1}{\sqrt{\lamone}} \, \left[ \ {\rm \|\rho\|_2 \pls \| \omega \|_2  \ } \right] \pls \frac{1}{\sqrt{\delone}} \, 
\left[ \ {\rm  \| \etanu \|_{2,\bdy}  \pls C_0 \| \omega \|_2 \  } \right]. \eeq
\etm 
\bpf
This proof is essentially the same as that of theorem \ref{T6.1}. 
The compatibility condition \eqref{e7.1} implies that compatibility condition for the solvability of \eqref{e5.9} with $\eta$
given by the right hand side of \eqref{e7.3} holds. 
The estimates now follow as in theorem \ref{T6.1}.
\epf

In a similar manner to corollary \ref{C6.3} of the last section one has   
\begin{cor} \label{C7.2}
\quad Suppose (B1) holds and $\curl \, \vbf, \Div \, \vbf \, {\rm \in \Lt,} \ \vbf {\rm \cdot \tau \in \Ltby}$, then $\vbf \in {\rm \Hdiv \cap \Hbcurl}$ 
and there is a  $C > 0$ such that
\beq \label{e7.8}
\| \vbf \|_2^2 \leqs {\rm C \, \left[ \ \|\, \curl \, \vbf \|_2^2 \pls \| \Div \, \vbf \,  \|_2^2   \pls   \| \vbf \cdot \tau \|_{2,\bdy}^2  \ \right]. } \eeq 
\end{cor}

When the region $\Om$ has holes, (that is its boundary has more than one connected component,  then the solution of this boundary value problem is non-unique. 
There are non-zero harmonic vector fields associated with potential differences between different components of the boundary. 
This was studied in \cite{AA} where the well-posed problem was described and  the solution described in the above theorem -  is the
least energy ( 2-norm) solution of the problem.
% In this case the norm of the harmonic fields is proportional to the relative capacity of the holes in the region. 

\vspace{2em}
%%%%%
%%%%%   line 1320
%%%%%

\section{Mixed boundary conditions}   \label{MBC}

In electromagnetic field theory, problems where given flux conditions are prescribed on part of the boundary
and tangential boundary data is prescribed on the complementary part need to be solved. 
The well-posedness and uniqueness of solutions of these problems was studied in sections 12 - 15 of Alexander and Auchmuty
\cite{AA}.
Here our primary interest is in obtaining 2-norm bounds on solutions in terms of the data. 
The constant in the relevant estimate will be the value of a natural optimization problem that is related to an eigenvalue in the 
case where the data is $L^2$.

The analysis of such problem differs considerably from that for the normal and tangential boundary value problems. 
First   no compatibility conditions on the data are required for $L^2-$solvability .
In addition the two potentials are each found directly by solving similar variational problems on  appropriate closed subspaces 
of $\Hone$ that is determined by the topology of the boundary data. 
The aim here  is to obtain bounds on the solutions of these problems, in particular of their energy,  in terms of norms of the data.
% A particular interest is what quantities enter into these bounds. 

The mixed div-curl boundary value problem is to find  vector fields $\vbf \,  {\rm \in \Ltto}$ that satisfy
\eqref{eDC.1} subject to
\beq   \label{BdM}
\vbf{\it (z) \cdot \tau(z)}   \eqs  {\rm \etatau(z)     \xon \Gamtau }\xand \vbf{\it (z) \cdot \nu(z)}   \eqs {\rm  \etanu(z)     \xon \Gamnu. }  \eeq
Here $\Gamtau, \Gamnu$ are nonempty open subsets of $\bdy$ whose union is dense in $\bdy$.
Our analysis will use the following requirement on these sets. 

\noindent{\bf Condition B2. }  {\it $\Gamma$ is an nonempty open subset of $\bdy$  with a finite number of disjoint components
$\{ \gamo, \ldots, \gamL \}$ and there is a finite distance $d_0$ such that $d(\gamj, \gamk) \geq d_0$ when $j \neq k$.}

When $\Gamma$ satisfies (B2), define $\HGamz$ to be the subspace of $\Hone$ of functions whose traces are zero on the
set $\Gamma \subset \bdy$. 
This is a closed subspace of $\Hone$ from lemma 12.1 of \cite{AA}.
When (B1) and (B2) hold then $\HGamz$ is a real Hilbert space with the $\pal-$inner product \eqref{ip2} 
Note that this inner product reduces to 
\beq \label{e8.3}
\ang{\vap, \psi}_{\pal, \Gamma} \deqs  \IOm \  \gradphi  \cdot \nabla \psi \  \dtx \pls  \IGamt \  \vap \; \psi \ ds \eeq
where $\Gamtil$ is the complement of $\Gamma$ in $\bdy$.

\begin{lem} \label{L8.1}
\quad Suppose $\Om, \bdy, \Gamma$ satisfy (B1) and (B2). Then \\
(i) \quad $\ang{\vap, \psi}_1 \deqs \ang{\gradphi, \gradpsi}$ is an equivalent inner product on $\HGamz$ to the $\pal-$inner product. \\
(ii) \quad For $q \in (1,\infty)$, there is an $M_q(\Gamma) $ such that 
\beq \label{e8.5}
\| \vap \|_q^q \pls \| \vap \|_{q,\bdy}^q \leqs M_q(\Gamma) \ \| \gradphi \|_2^q \foral \vap \in \HGamz. \eeq
(iii) \quad if $\Gamma_1 \supset \Gamma$, then $M_q(\Gamma_1)  \leqs M_q(\Gamma)$.  
\end{lem}   \bpf
(i) \quad Let $\lamogam$ be the least eigenvalue of the Laplacian on $\HGamz$ so that 
\beq \label{e8.7}
\IOm \ | \gradphi |^2 \dtx \geqs \lamogam \ \IOm \ \vap^2 \, \dtx \foral \vap \in \HGamz. \eeq
This  exists and is positive when (B2) holds as $\sigma(\Gamma) > 0$;
see proposition 13.2 in \cite{AA} for a proof.
Since the $\pal-$ norm and the standard norm on $\Hone$ are equivalent, there is a $C > 0$ such that 
\[ \| \vap \|_{\pal}^2 \leqs C \, \| \vap \|_{1,2}^2 \leqs C \, \left( \, 1 \pls \lamogam^{-1} \, \right) \, \| \gradphi \|_2^2 \]
Thus the norm from (i) is equivalent to the $\pal$-norm. 

(ii) \quad Consider the functional $\Gcq (\vap) \deqs \| \vap \|_q^q \pls \| \vap \|_{q,\bdy}^q $ on $\HGamz$.
This functional is convex and weakly continuous as the imbedding of $\Hone$ into $\Lq$ and $\Lqby$ are compact for any $q \geq 1$
when (B1) holds. 
Let $B_1$ be the unit ball in $\HGamz \ \wrt$ the inner product of (i). 
Define $M_q(\Gamma) := \sup_{\vap \in B_1} \ \Gcq(\vap)$.
This sup is finite and \eqref{e8.5} follows upon scaling.

(iii) \quad When $\Gamma_1 \supset \Gamma$, then $H^1_{\Gamma_10} (\Om) \subset \HGamz$, so the associated unit ball
is smaller and thus $M_q(\Gamma_1)  \leqs M_q(\Gamma)$.
\epf

It appears that the value of $M_q(\Gamma)$ increases to $\infty$ as $\sg(\Gamma)$ decreases to zero. 
It would be of interest  to estimate or quantify  this dependence.
When $q=2$ the constant $M_2(\Gamma)$ is related to the least eigenvalue  of an eigenvalue problem
for the Laplacian where the eigenvalue appears in both the equation and the boundary condition. 
Suppose that $\lambda_1(\Om, \Gamma)$ is the least eigenvalue of
\beq \label{e8.8}
- \Delta u \eqs \lambda \, u \xon \Om \xwith u = 0 \xon \Gamma, \quad  \Dnu u = \lambda u \xon \Gamtil, \eeq
then $M_2(\Gamma) \eqs  \lambda_1(\Om, \Gamma)^{-1}.$

 The conditions required here are that $\Gamtau, \Gamnu$ are proper subsets of $\bdy$ satisfying the following. 

 \noindent{\bf Condition B3. }  {\it $\Gamnu$ and $\Gamtau$ are disjoint, satisfy (B2) and have  union  dense in $\bdy$.}

Let  $\GradGam, \CuGam$  be the spaces of gradients of functions in $\HGamtau$ and curls of functions in 
$\HGamnu$ respectively. 
These spaces are $L^2-$ orthogonal as fields in $\Ltto$.
The vector field $\vbf := \gradp \psi - \gradphi$ will satisfy the boundary condition \eqref{BdM} in a weak sense provided
$ \vap \in \HGamnu$ and $\Dnu \vap  + \etanu = 0$ on $\Gamnu$ and $\psi \in \HGamtau$ with $\Dnu \psi + \etatau = 0$ 
on $\Gamtau$.

As described in \cite{AA} there are variational principles for the potentials in this representation and the  field $\vbf$ will be a solution of 
\eqref{eDC.1} - \eqref{BdM} provided $\vap \in \HGamtau$ is a solution of 
\beq \label{e8.10}
\IOm [ \, \gradphi \cdot \gradchi \mns \rho \chi \, ] \, \dtx \pls \IGamnu \etanu \chi \, \dsg \eqs 0 \foral \chi \in \HGamtau. \eeq 
Similarly  $\psi \in \HGamnu$ is a solution of 
\beq \label{e8.11}
\IOm [ \, \gradpsi \cdot \gradchi \mns \omega \chi \, ] \, \dtx \pls \IGamtau \etatau \chi \, \dsg \eqs 0 \foral \chi \in \HGamnu. \eeq 

Note that these equations are of the same type; they differ only in that $\Gamtau, \Gamnu$ are interchanged from one to the other.
They can be written as a problem of finding  $\vap \in \HGamz$ satisfying
\beq \label{e8.12} 
\IOm \ \gradphi \cdot \gradchi \ \dtx \eqs \Fc (\chi) \foral \chi \in \HGamz. \eeq
Here $\Fc(\chi) $ is the linear functional defined by $\Fc(\chi) \eqs \IOm \rho \, \chi \, \dtx \mns \Iby \eta \, \chi \,  ds$.
For notational convenience  the functions $\etanu, \etatau$ are extended to all of  $\bdy $ by zero.

The general result about this problem may be described as follows.
\btm \label{T8.2}
Assume that $\Om, \Gamma$ satisfy (B1)-(B2) with $\rho \in  \Lq, \eta  \in \Lqby, q > 1$.
Then there is a unique solution $\phit \in \HGamz$ of \eqref{e8.12} and it satisfies
\beq \label{e8.15}
\| \nabla \phit \|^q_2  \leqs   M_{q'}(\Gamma)^{q-1} \, \left[ \ \|\, \rho \, \|^q_q \pls \|\, \eta \, \|^q_{q, \Gamtil} \right] \eeq
\etm
\bpf
When $\phit$ is a solution of \eqref{e8.12} and $|\Fc(\chi) | \leq C \, \| \gradchi \|_2 \ $ for all $\ \chi \in \HGamz$, then 
$\| \nabla \phit \|_2 \leq C$. 
So the result just requires an appropriate estimate of $\Fc(\chi)$.
Two applications of Holder's inequality to the definition of $\Fc$  yield that 
\[ |\Fc(\chi)| \leqs \left[ \, \|\rho\|_q^q \pls  \|\eta\|_{q,\bdy}^q \, \right]^{1/q}  . \left[ \, \|\chi\|_{q'}^{q'} \pls  \|\chi \|_{q',\bdy}^{q'} \, \right]^{1/q'} \]
for all $\chi \in \HGamz$.
Then \eqref{e8.3} yields 
\[ |\Fc(\chi)| \leqs \left[ \, \|\rho\|_q^q \pls  \|\eta\|_{q,\bdy}^q \, \right]^{1/q}  M_{q'}(\Gamma)^{1/q'} \  \| \gradchi \|_2 . \]
This inequality yields \eqref{e8.15}.
\epf

\begin{cor}
Assume that $\Om, \Gamnu, \Gamtau$ satisfy (B1)-(B3) with $\rho, \omega \in \Lq, \etanu, \etatau \in \Lqby$ and $ q > 1$.
Then there is a solution $\tilde{\vbf} = \gradp \psit  \mns \nabla \phit$ of \eqref{eDC.1} - \eqref{BdM} with
\beq \label{e9.15}
\| \, \tilde{\vbf} \, \|^2_2  \  \leqs C_q(\Gamtau) \left[ \| \rho \|^q_q \pls  \| \etanu \|^q_{q, \Gamnu} \right] ^{2/q} \pls 
C_q(\Gamnu) \left[ \| \omega \|^2_q \pls  \| \etatau \|^2_{q, \Gamtau} \right]^{2/q}     \eeq
If $\vbf$ is any solution of this mixed div-curl system, then $\| \vbf \|_2 \geqs \| \, \tilde{\vbf} \, \|_2$.
\end{cor}
\bpf
Let $\phit, \psit$ be the solutions of \eqref{e8.10} - \eqref{e8.11} respectively.
Then their orthogonality implies that $\| \, \tilde{\vbf} \, \|^2_2 \eqs \| \nabla \phit \|^2_2 \pls \| \nabla \psit \|^2_2 $.
Theorem \ref{T8.2} implies that there is a constant such that 
\[  \| \nabla \phit \|^2_2 \leqs C_q(\Gamtau) \ \left[ \ \|\, \rho \, \|^q_q \pls \|\, \etanu \, \|^q_{q, \Gamnu} \right]^{2/q}. \] 

Similarly the other mixed boundary value problem has solution $\psit$ with
\[  \| \nabla \psit \|^2_2 \leqs C_q(\Gamnu) \ \left[ \ \|\, \omega \, \|^q_q \pls \|\, \etatau \, \|^q_{q, \Gamtau} \right]^{2/q}. \] 
Adding these two expressions leads to the inequality of \eqref{e9.15}.

\epf

In general there is an affine subspace of solutions of \eqref{eDC.1} - \eqref{BdM}  as described in section 14 of \cite{AA}.
To find a well-posed problem certain linear functionals   of the solutions must be further specified and the energy of the
solution depends on these extra imposed conditions.

%\vspace{3em}
%%%%%   line 1500
%%%%%

% \newpage
%
% ########    Bibliography  #########
%
\vspace{2em}

Since these are  problems of interest to researchers in a variety of different areas, some relevant references to the literature may have been omitted from the following bibliography.
 The author would appreciated being informed about  further papers that treat these  topics  analytically.

%\vspace {5em}

\end{document}